\documentstyle[12pt]{article}
\newtheorem{theorem}{Theorem}\newtheorem{lemma}{Lemma}
\newtheorem{corollary}{Corollary}

\newtheorem{conjecture}{Conjecture}
\newtheorem{definition}{Definition}

\def\pont{\hspace{-6pt}{\bf.\ }}   
\def\beq{\begin{equation}}\def\eeq{\end{equation}}
\def\beqn{\begin{eqnarray}}\def\eeqn{\end{eqnarray}}

\def\eps{\varepsilon}
\def\({\mbox{$($}}\def\){\mbox{$)$}}

\def\qed{\ifhmode\unskip\nobreak\fi\quad\ifmmode\Box\else$\Box$\fi}
\newcommand{\text}[1]{\quad\mbox{#1}\quad}

\title{ Monochromatic bounded degree subgraph partitions }

\author{Andrey Grinshpun
\thanks{Research supported in part by
National Physical Sciences Consortium Fellowship}\\[-0.8ex]
\small Department of Mathematics\\[-0.8ex]
\small Massachusetts Institute of Technology\\[-0.8ex]
\small Cambridge, MA, USA 02139\\[-0.8ex]
\small \texttt{agrinshp@mit.edu}
 \and G\'{a}bor N. S\'ark\"ozy
\thanks{Research supported in part by
OTKA Grant No. K104343.}\\[-0.8ex]
\small Alfr\'ed R\'enyi Institute of Mathematics\\[-0.8ex]
\small Hungarian Academy of Sciences\\[-0.8ex]
\small Budapest, P.O. Box 127\\[-0.8ex]
\small Budapest, Hungary, H-1364\\
and\\
\small Computer Science Department\\[-0.8ex]
\small Worcester Polytechnic Institute\\[-0.8ex]
\small Worcester, MA, USA 01609\\[-0.8ex]
\small \texttt{gsarkozy@cs.wpi.edu}}

\begin{document}
\maketitle

\begin{abstract}
Let ${\cal{F}}=\{F_1,F_2,\ldots\}$ be a sequence of graphs such that
$F_n$ is a graph on $n$ vertices with maximum degree at most
$\Delta$. We show that there exists an absolute constant $C$ such
that the vertices of any 2-edge-colored complete graph can be
partitioned into at most $2^{C\Delta \log{\Delta}}$ vertex disjoint
monochromatic copies of graphs from ${\cal{F}}$. If each $F_n$ is
bipartite, then we can improve this bound to $2^{C \Delta}$; this
result is optimal up to the constant $C$.
\end{abstract}

\section{Introduction}

Let $K_n$ be a complete graph on $n$ vertices whose
edges are colored with $r$ colors ($r\geq 1$). How many
monochromatic cycles (single vertices and edges are considered to be cycles) are needed to partition the vertex set of
$K_n$? This question received much attention in the last few
years. Let $p(r)$ denote the minimum number of
monochromatic cycles needed to partition the vertex set of any
$r$-colored $K_n$. It is not obvious that $p(r)$ is a well-defined
function. That is, it is not obvious that there always is a
partition whose cardinality is independent of $n$. However, in \cite{EGYP} Erd\H{o}s, Gy\'arf\'as, and
Pyber proved that there exists a constant $C$ such that $p(r) \leq C
r^2 \log r$ (throughout this paper $\log$ denotes the natural
logarithm). Furthermore, in \cite{EGYP} (see also \cite{GY}), the
authors conjectured that $p(r)=r$.

The special case $r=2$ of this conjecture was asked earlier by Lehel
and for $n\geq n_0$ was first proved by {\L}uczak, R\"{o}dl, and
Szemer\'edi \cite{LRS}. Allen improved on the value of $n_0$
\cite{A} and recently Bessy and Thomass\'e \cite{BT} proved the
original conjecture for $r=2$. For general $r$ the current best
bound is due to Gy\'arf\'as, Ruszink\'o, S\'ark\"ozy, and
Szemer\'edi \cite{GRSS} who proved that for $n\geq n_0(r)$ we have
$p(r) \leq 100r\log r$. For $r=3$ an approximate version of the
conjecture was proved in \cite{GYSSZ2} but, surprisingly, Pokrovskiy
\cite{PO} found a counterexample to the conjecture. However, in the
counterexample, all but one vertex can be covered by $r$ vertex
disjoint monochromatic cycles. Thus, a slightly weaker version of
the conjecture still can be true, say that, apart from a constant
number of vertices, the vertex set can be covered by $r$ vertex
disjoint monochromatic cycles.

Let us also note that the above problem was generalized in various
directions; for hypergraphs (see \cite{GS} and \cite{SG}), for
complete bipartite graphs (see \cite{EGYP} and \cite{H}), for graphs
which are not necessarily complete (see \cite{BBGGS} and \cite{SA}),
and for partitions by monochromatic connected $k$-regular graphs
(see \cite{SS} and \cite{SSS}).

Another area that attracted a lot of interest is the study of Ramsey
numbers for bounded degree graphs. For a graph $G$, the {\em Ramsey
number} $R(G)$ is the smallest positive integer $N$ such that if the
edges of a complete graph $K_N$ are partitioned into two color
classes then one color class has a subgraph isomorphic to $G$. The
existence of such a positive integer is guaranteed by Ramsey's
classical result \cite{RA}. Determining $R(G)$ even for very special
graphs is notoriously hard (see e.g. \cite{GRS} or \cite{RAD}).

In 1975, Burr and Erd\H{o}s \cite{BE} raised the problem that every
graph $G$ with $n$ vertices and maximum degree $\Delta$ has a linear
Ramsey number, so $R(G) \leq C(\Delta)n$, for some constant
$C(\Delta)$ depending only on $\Delta$. This was proved by
Chv\'atal, R\"odl, Szemer\'edi and Trotter \cite{CRST} in one of the
earliest applications of Szemer\'edi's celebrated Regularity Lemma
\cite{Sz}. Because the proof uses the Regularity Lemma, the bound on
$C(\Delta)$ is quite weak; it is of tower type in $\Delta$. This was
improved by Eaton \cite{E}, who proved, using a variant of the
Regularity Lemma, that the function $C(\Delta)$ can be taken to be
of the form $2^{2^{O(\Delta)}}$.

Soon after, Graham, R\"odl, and Ruci\'nski \cite{GRR} improved this
further to $C(\Delta)\leq 2^{O(\Delta\log^2{\Delta})}$ and for
bipartite graphs $C(\Delta)\leq 2^{O(\Delta\log{\Delta})}$. They
also proved that there are bipartite graphs with $n$ vertices and
maximum degree $\Delta$ for which the Ramsey number is at least
$2^{\Omega(\Delta)}n$. Recently, Conlon \cite{C} and, independently,
Fox and Sudakov \cite{FS} have shown how to remove the
$\log{\Delta}$ factor in the exponent, achieving an essentially best
possible bound of $R(G) \leq 2^{O(\Delta)} n$ in the bipartite case.
For  the non-bipartite graph case, the current best bound is due to
Conlon, Fox, and Sudakov \cite{CFS} $C(\Delta)\leq
2^{O(\Delta\log{\Delta})}$. Similar results have been proven for
hypergraphs: \cite{CFKO1,CFKO2,NORS} use the Hypergraph Regularity
Lemma and \cite{CFS1} improves the bounds by avoiding the Regularity
Lemma.

It is a natural question (initiated by Andr\'as Gy\'arf\'as) to
combine the above two problems and ask how many monochromatic
members from a bounded degree graph family are needed to partition
the vertex set of a 2-edge-colored $K_N$. In this paper we study
this problem. Given a sequence ${\cal{F}}=\{F_1,F_2,\ldots\}$ of
graphs, we say it is \emph{$\Delta$-bounded} if each $F_n$ is a
graph on $n$ vertices with maximum degree at most $\Delta$. In
general we say that ${\cal{F}}$ has some graph property if every
graph of $\cal F$ has that property (e.g. ${\cal F}$ is bipartite if
$F_n$ is bipartite for every $n$).

We prove the following result on partitions by monochromatic members
of ${\cal{F}}$.
\begin{theorem}\pont\label{tetel}
There exists an absolute constant $C$ such that, for every $\Delta$
and every $\Delta$-bounded graph sequence $\cal{F}$, every
$2$-edge-colored complete graph can be partitioned into at most
$2^{C\Delta\log{\Delta}}$ vertex disjoint monochromatic graphs from
$\cal{F}$.
\end{theorem}

Thus, perhaps surprisingly, we have the same phenomenon as for
cycles; we can partition into monochromatic graphs from ${\cal{F}}$
such that their average size is roughly the same as the single
largest monochromatic graph we can find. In the case of a bipartite
${\cal{F}}$ we can eliminate the $\log{\Delta}$ factor from the
exponent to get the following essentially best possible result.

\begin{theorem}\pont\label{bip-tetel}
There exists an absolute constant $C$ such that, for every $\Delta$
and every bipartite $\Delta$-bounded graph sequence $\cal F$, every
$2$-edge-colored complete graph can be partitioned into at most
$2^{C\Delta}$ vertex-disjoint monochromatic copies of graphs from
$\cal{F}$.
\end{theorem}

We do not make an effort to optimize the constant $C$ since probably
it will be far from optimal anyway. However, in both theorems we
must use at least $2^{\Omega(\Delta)}$ parts.

\begin{theorem}\pont\label{lower}
There exists an absolute constant $c$ such that, for every $\Delta$,
there is a bipartite $\Delta$-bounded graph sequence $\cal F$ and
there is a $2$-edge-coloring of $K_n$ so that covering the vertices
of $K_n$ using monochromatic copies of graphs from $\cal F$ requires
at least $2^{c\Delta}$ such copies.
\end{theorem}

It would be desirable to close the gap between the upper and lower
bounds for non-bipartite ${\cal{F}}$ as well, though doing so may
require improved bounds for the Ramsey numbers of bounded degree
graphs. Furthermore, it would be interesting to extend this problem
for more than 2 colors.

Let us also mention one interesting special case of our theorem. The
$k^{th}$ {\em power} of a cycle $C$ is the graph obtained from $C$
by joining every pair of vertices with distance at most $k$
in $C$. Density questions for powers of cycles
have generated a lot of interest; in particular the famous
P\'osa-Seymour conjecture (see e.g.
\cite{CDK,FH,FK2,FK3,FGJS,KSSz1,KSSz4,KSSz5,LSSz}). Theorem
\ref{tetel} implies the following result on the partition number by
monochromatic powers of cycles.

\begin{corollary}\pont
There exists an absolute constant $C$ so that for every $k$ every
$2$-colored complete graph can be partitioned into at most
$2^{Ck\log{k}}$ vertex disjoint monochromatic $k$th powers of
cycles.
\end{corollary}

However, we must note that in this case probably the optimal answer
is $O(k)$.

\section{Notation and tools}

For basic graph concepts see the monograph of Bollob\'as \cite{B}.\\
$V(G)$ and $E(G)$ denote the vertex-set and the edge-set of the
graph $G$. $(A,B,E)$ denotes a bipartite graph $G=(V,E)$, where
$V=A\cup B$ and $E\subset A\times B$. A proper $r$-coloring of $G$ is a coloring of its vertices where no two adjacent vertices receive the same color. For a graph $G$ and a subset
$U$ of its vertices, $G|_U$ is the restriction to $U$ of $G$. $N(v)$
is the set of neighbors of $v\in V$. Hence, $|N(v)|=deg(v)=deg_G(v)$,
the degree of $v$. $\delta(G)$ stands for the minimum and
$\Delta(G)$ for the maximum degree in $G$. When $A,B$ are subsets of
$V(G)$, we denote by $e(A,B)$ the number of edges of $G$ with one
endpoint in $A$ and the other in $B$. In particular, we write
$deg(v,U)=e(\{v\},U)$ for the number of edges from $v$ to $U$. For
non-empty $A$ and $B$,$$d(A,B)=\frac{e(A,B)}{|A||B|}$$ is the {\em
density} of the graph between $A$ and $B$.
\begin{definition}\pont
The bipartite graph $G=(A,B,E)$ is $\eps$-regular if
$$X\subset A,\ Y\subset B,\ |X|>\eps|A|,\ |Y|>\eps|B|
\text{imply}|d(X,Y)-d(A,B)|<\eps.$$
\end{definition}
We will often say simply that ``the pair $(A,B)$ is $\eps$-regular''
with the graph $G$ implicit.
\begin{definition}\pont $(A,B)$ is
$(\eps,d,\delta)$-super-regular if
it is $\eps$-regular, satisfies $d(A,B)\geq d$, and
$$deg(a)\geq \delta|B|\ \forall \; a\in A,\qquad
deg(b)\geq \delta|A|\ \forall \; b\in B.$$ An
$(\eps,\delta,\delta)$-super-regular pair is simply called
$(\eps,\delta)$-super-regular (and in this case the density
condition is not needed).
\end{definition}

We will use frequently the following well-known property of regular
pairs claiming that subsets of a regular pair also form a regular
pair with somewhat weaker parameters.
\begin{lemma}[Slicing Lemma, Fact 1.5 in \cite{KS}]\pont\label{regu}
Let $(A,B)$ be an $(\eps, d,0)$-super-regular pair (i.e. one with no
minimum degree constraint), and, for some $\beta > \eps$, let
$A'\subset A$, $|A'|\geq \beta |A|$, $B'\subset B$, $|B'|\geq \beta
|B|$. Then $(A',B')$ is an $(\eps', d',0)$-super-regular pair with
$\eps'=\max \{\eps/\beta, 2\eps \}$ and $|d' - d|< \eps$.
\end{lemma}

\begin{definition}\pont Given a $k$-partite graph $G=(V,E)$ with $k$-partition $V=V_1\cup \ldots \cup V_k$, the $k$-cylinder
$V_1\times \ldots \times V_k$ is $\eps$-regular
($(\eps,d,\delta)$-super-regular or $(\eps,\delta)$-super-regular)
if all the ${k \choose 2}$ pairs of subsets $(V_i,V_j)$, $1\leq i <
j\leq k$, are $\eps$-regular ($(\eps,d,\delta)$-super-regular or
$(\eps,\delta)$-super-regular). Given $0< \alpha <1$, the
$k$-cylinder $V_1\times \ldots \times V_k$ is $\alpha$-balanced if,
for every $i < j$, $||V_i|-|V_j|| \leq \alpha min(|V_i|,|V_j|)$.
\end{definition}

Instead of the Regularity Lemma of Szemer\'edi \cite{Sz} we will use
the following lemma which Conlon and Fox \cite{CF} argued as a
consequence of the Duke, Lefmann, and R\"odl weak Regularity Lemma
\cite{DLR}.

\begin{lemma}[\cite{DLR} and Lemma 5.3 in \cite{CF}]\pont\label{fox}
For each $0 < \eps < 1/2$, any graph $G = (V,E)$ on at least $k$
vertices has an $\eps$-regular $k$-cylinder with parts of equal size
(i.e. 0-balanced); the size of each part is at least $\frac{1}{2k}
\eps^{k^2\eps^{-5}} |V|.$
\end{lemma}

We will use the following corollary of this lemma.

\begin{lemma}[Lemma 5.4 in \cite {CF}]\pont\label{reg}
For each $0 < \eps < 1/2$, any 2-colored complete graph on at least
$2^{2k}$ vertices has, in one of the colors (say in red), an
$(\eps,1/2,0)$-super-regular $0$-balanced $k$-cylinder (i.e. one
with no minimum degree constraint and parts of equal size), where
the size of each part is at least $\frac{1}{2(2^{2k})}
\eps^{2^{4k}\eps^{-5}} n$.
\end{lemma}

Indeed, to get this one applies Lemma \ref{fox} for the red subgraph
with $2^{2k}$ in place of $k$ to get an $\eps$-regular
$2^{2k}$-cylinder. Then we may consider the complete graph whose
vertices $i$ correspond to the parts of the cylinder $V_i$ and we
color the edge $(i,j)$ by the majority color in the pair
$(V_i,V_j)$. We then apply $R(K_k)\leq 2^{2k}$ and use the fact
that, if $(V_i,V_j)$ is regular in  one color, then it is also
regular in the other color.

We will also use the Hajnal-Szemer\'edi Theorem on equitable proper
colorings. A proper coloring is {\em equitable} if the numbers of vertices
in any two color classes differ by at most one.

\begin{lemma}[Hajnal-Szemer\'edi Theorem \cite{HSZ}]\pont\label{h-sz}
Any graph with maximum degree at most $\Delta$ has an equitable proper
coloring with $\Delta + 1$ colors.
\end{lemma}
For a simpler proof see also \cite{KK}.

Our other main tool is a quantitative version of the Blow-up Lemma
(see \cite{KSSz2,KSSz3}).
\begin{lemma}[Quantitative Blow-up Lemma]\label{blowup}\pont
There exists an absolute constant $C_{BL}$ such that, given a graph
$R$ of order $r\geq 2$ and positive parameters $\delta$ and
$\Delta$, for any $0< \eps < (\frac{\delta
d^{\Delta}}{r\Delta})^{C_{BL}}$ the following holds. Let $N$ be an
arbitrary positive integer, and let us replace the vertices of $R$
with pairwise disjoint $N$-sets $V_1,V_2,\ldots,V_r$ (blowing up).
We construct two graphs on the same vertex-set $V=\bigcup V_i$. The
graph $R(N)$ is obtained by replacing all edges of $R$ with copies
of the complete bipartite graph $K_{N,N}$, and a sparser graph $G$
is constructed by replacing the edges of $R$ with some
$(\eps,d,\delta)$-super-regular pairs. If a graph $H$ with
$\Delta(H)\leq\Delta$ is embeddable into $R(N)$, then it is
embeddable into $G$.
\end{lemma}

Thus, roughly speaking, regular cylinders behave as complete partite
graphs from the viewpoint of embedding bounded degree subgraphs.
Note that, in either of the proofs \cite{KSSz2, KSSz3}, the
dependence of $\eps$ on the other parameters was not computed
explicitly. To prove Lemma \ref{blowup} one has to go through the
proof, say from \cite{KSSz3}, and make all the dependencies of the
parameters explicit. All the details are presented in \cite{QBL}.
Note that for the proof of Theorem \ref{tetel} we could use
$d=\delta$, we need the stronger version for the proof of Theorem
\ref{bip-tetel}, as $\delta$ may be much smaller than $d$.

In particular we will need the following consequence of the Blow-up
Lemma.

\begin{lemma}\label{cylinder}\pont
There exists an absolute constant $C_{BL}$ such that, given positive
parameters $\delta$, $d$, and $\Delta$, and given a $\Delta$-bounded graph
sequence ${\cal{F}}$, for any $0< \eps <
\frac{1}{2}\left(\frac{(\frac{\delta}{2})
(\frac{d}{2})^{\Delta}}{\Delta(\Delta+2)}\right)^{C_{BL}}$ the
following holds. Let $G=(V,E)$ be a $(\Delta+2)$-partite graph with
$(\Delta+2)$-partition $V=V_1\cup \ldots \cup V_{\Delta+2}$, where
the cylinder $V_1\times \ldots \times V_{\Delta+2}$ is
$(\eps,d,\delta)$-super-regular and $\eps$-balanced. Then we can
partition the vertex set into at most $(\Delta+3)$ vertex disjoint
copies of graphs from ${\cal{F}}$.
\end{lemma}

Indeed, note first that by the Blow-up Lemma (Lemma \ref{blowup}
with $r=\Delta + 2$), if the cylinder is $0$-balanced, then it is
enough to check the statement for the complete $(\Delta +2)$-partite
graph with the same partite sets. But then the Hajnal-Szemer\'edi
Theorem (Lemma \ref{h-sz}) implies that we may cover the cylinder with a
single graph from ${\cal{F}}$ (even if we had only $\Delta + 1$
partite sets instead of $\Delta + 2$).

If the cylinder is not $0$-balanced (but $\eps$-balanced), then
first we eliminate the small discrepancies among the sizes. For each
index $i \in [\Delta + 2]$ define $v_i = |V_i|$ and take $v =
\max_i(v_i)$. Then, for each set $S \subseteq [\Delta + 2]$ of size
$\Delta + 1$, define $w_S = v-v_i$ where $i$ is the unique index not
contained in $S$. For each such set $S$, by the Blow-up Lemma, we
may find a copy of a graph from $\cal{F}$ that uses $w_S$ vertices
from each $V_i$ with $i \in S$. Thus, we use $(\Delta + 2)$ such
graphs, one for each $S$.

After this procedure, the number of vertices remaining in $V_i$ is
\beq\label{v}v_i - \sum_{S : i \in S} w_S = v_i - \sum_S w_S +
w_{[\Delta + 2] \setminus \{i\}} = v - \sum_S w_S \;\; \left(\geq
(1- (\Delta + 2) \eps) v\right).\eeq That is, after this procedure,
the cylinder is $0$-balanced, and we may cover the remaining
vertices with a single graph from $\cal{F}$ by the Blow-up Lemma.
Indeed, the Slicing Lemma (Lemma \ref{regu}) and (\ref{v}) imply
that the remaining cylinder is $(2\eps,d/2,\delta/2)$-super-regular
and $0$-balanced so we may indeed apply the Blow-up Lemma (Lemma
\ref{blowup}). We use a total of at most $(\Delta + 3)$ graphs.

\section{Proof of Theorem \ref{tetel}}

For technical reasons, it will be convenient to prove a stronger
version of Theorem \ref{tetel}.

\begin{theorem}\pont\label{tetel2}
There exists an absolute constant $C$ so that if ${\cal{F}}_1$ and
${\cal{F}}_2$ are $\Delta$-bounded graph sequences and every graph
in ${\cal{F}}_1, {\cal{F}}_2$ has an equitable proper coloring with
$\chi_1,\chi_2$ colors, respectively, then every $2$-edge-colored
complete graph can be partitioned into at most $2^{C(\chi_1 + \chi_2
+ \Delta)\log{\Delta}}$ vertex disjoint red copies of
graphs from ${\cal{F}}_1$ and monochromatic blue copies of graphs
from ${\cal{F}}_2$.
\end{theorem}

To see that this implies Theorem \ref{tetel}, note again that by the
Hajnal-Szemer\'edi Theorem (Lemma \ref{h-sz}) any graph of maximum
degree at most $\Delta$ has an equitable proper coloring with
$\Delta + 1$ colors.

To avoid redundancy, for the rest of this section ${\cal{F}}_1$ and
${\cal{F}}_2$ will be $\Delta$-bounded sequences.

We have three main tools for finding monochromatic copies of graphs.
One of these is the Blow-up Lemma (Lemma \ref{cylinder}). Another is
applying bounds on Ramsey numbers:

\begin{lemma}\label{ramsey}\pont
There exists an absolute constant $C_1$ such that given a 2-edge-colored
$K_n$ and an $\eps > 0$, we may cover all but an $\eps n$ fraction
of the vertices using vertex disjoint monochromatic red copies of
graphs from ${\cal{F}}_1$ and blue copies of graphs from
${\cal{F}}_2$ while using at most $2^{C_1\Delta
\log{\Delta}}\log(1/\eps)$ such copies.
\end{lemma}

\noindent {\bf Proof. } By the bound in \cite{CFS} (see Theorem
2.6), any 2-edge-colored $K_n$ contains either a red copy of $F_1$ or a
blue copy of $F_2$ for any two graphs $F_1,F_2$ with maximum degree
at most $\Delta$ that are on $n'=2^{-C_1\Delta \log{\Delta}}n$
(assume for simplicity that this is an integer) vertices. Pick any
graph in $F_1 \in {\cal{F}}_1$ and $F_2 \in {\cal{F}}_2$ with $n'$
vertices, find a copy in an appropriate color of one of these in our
$K_n$, remove the vertices of this copy from the graph, and recurse
on the remaining $n-n'$ vertices. Note the number of remaining
vertices is $(1-2^{-C_1\Delta \log{\Delta}})n \leq e^{-2^{-C_1\Delta
\log{\Delta}}}n$, so after repeating this process $2^{C_1\Delta
\log{\Delta}}\log(1/\eps)$ times, we are left with a graph on at
most $e^{-\log(1/\eps)}n = \eps n$ vertices, as desired. \qed

The final tool we have is a simple greedy embedding of bounded
degree bipartite graphs into very dense bipartite graphs.

\begin{lemma}\label{bipartite}\pont
Given a bipartite graph $H=(A,B,E)$ where every vertex of $B$ has
degree at most $\Delta$, a bipartite graph $G=(A',B',E')$ where
every vertex of $A'$ has degree at least $(1-1/(2\Delta))|B'|$ and
where $|B'| \geq 2|B|$, and an injection $\phi: A \rightarrow A'$,
the function $\phi$ extends to an injective homomorphism from $H$ to
$G$.
\end{lemma}

\noindent {\bf Proof. } We will embed the vertices of $B=\{b_1,b_2,\ldots,b_k\}$ one at a time (in this order). Once we have embedded $b_1,\ldots,b_i$, we show how to embed $b=b_{i+1}$. Note that $b$ must be embedded in a way consistent with its neighbors; i.e., $b$ must be contained in the common neighborhood (in $G$) of $\phi(N_H(b))$. Since $b$ has degree at most $\Delta$ and all of the vertices in $A'$ have degree at least $(1-\frac{1}{2 \Delta})|B'|$, by a union bound the number of vertices in $B'$ consistent with the neighbors of $b$ is at least
$$|B'|-\Delta \frac{1}{2\Delta}|B'| = \frac{|B'|}{2} \geq |B|.$$
Since we have embedded only $i < |B|$ vertices so far, at least one of the above $|B|$ vertices has not yet had any vertex embedded to it; to this vertex we embed $b$.
When this procedure embeds $b_k$, the embedding is complete, and is an injective homomorphism that extends $\phi$ by construction. \qed

We will combine these three tools to prove Theorem \ref{tetel2}. We
basically follow the {\em greedy-absorbing} proof technique that
originated in \cite{EGYP} and is used in many papers in this area
(e.g. \cite{GRSS}, \cite{H}, \cite{SSS}). We establish the desired
bound in the following steps.

\begin{itemize}
\item Step 1: First we find a special object, a regular cylinder, that is dense
in a color (say red) as given by Lemma \ref{reg}.
\item Step 2: Then we remove the cylinder and by iterating Ramsey's
Theorem as it was done in Lemma \ref{ramsey} we cover most of the
remaining vertices with monochromatic copies of graphs from
${\cal{F}}_1,{\cal{F}}_2$.
\item Step 3: Finally we add the few vertices that are neither in the cylinder nor
covered by monochromatic copies to the cylinder; since there will be
few vertices added, this will not affect the regularity of the
cylinder and the cylinder will absorb these vertices. Indeed, if it
were the case that all of the vertices in the cylinder had large
enough degrees in red to all the other partite sets, then the
Blow-up Lemma (Lemma \ref{cylinder}) would allow us to cover all of
them with a few red subgraphs from ${\cal{F}}_1$. By regularity,
there may be few vertices that fail to meet this minimum degree
condition, and they must have large degree in blue to one of the
sets in the cylinder. Inductively, we will partition these remaining
vertices into either red copies of graphs in ${\cal{F}}_1$ or blue
copies of graphs obtainable by taking a graph $F_2$ from
${\cal{F}}_2$ and removing an equitable color class. Then we will
use the fact that these vertices have large degree in blue and Lemma
\ref{bipartite} to ``glue in'' the missing parts of the $F_2$ graphs
by using some vertices from the cylinder.
\end{itemize}

We now proceed with proving Theorem \ref{tetel2}. The proof is by
induction on $\chi_1+\chi_2$. If either $\chi_1$ or $\chi_2$ is $1$,
then the result is trivial, as every graph in the corresponding
collection is an independent set. Otherwise, assume both $\chi_1$
and $\chi_2$ are at least $2$. Let any $2$-edge-colored $K_n$ be
given. We may assume that $n$ is at least $2^{2(\Delta + 2)}$ since
otherwise we can cover it by isolated vertices. We may also assume
$\Delta\geq 2$. Let \beq\label{para}\eps = 2^{-C_2\Delta},
k=\Delta+2,\mbox{ and } \eta = \frac{1}{2(2^{2k})}
\left(\frac{\eps}{2}\right)^{2^{4k}(\frac{\eps}{2})^{-5}},\eeq with
some sufficiently large absolute constant $C_2$ (independent of
$\Delta$). Apply Lemma \ref{reg} with parameter $\eps/2$ to the
2-coloring to get a 0-balanced $k$-cylinder $V = V_1\cup \cdots \cup
V_k$, where, for each $i, 1 \leq i \leq k$, we have $|V_i| \geq \eta
|V|$, and there is a color, say red, so that the cylinder is
$(\eps/2,1/2,0)$-super-regular in the red subgraph.

By Lemma \ref{ramsey}, using at most $2^{C_1\Delta
\log{\Delta}}\log(2/(\eps^2 \eta))$ copies, we may cover all but
$\eps^2 \eta n/2$ vertices of $K_n\setminus V$ with monochromatic
graphs in the appropriate color from ${\cal{F}}_1,{\cal{F}}_2$. Take
the remaining $\eps^2 \eta n/2$ vertices and add them to the
cylinder in such a way that the cylinder remains as balanced as
possible. Note that since we added at most $\eps^2 \eta n/2$
vertices in this way, the resulting cylinder $V' = V_1' \cup \cdots
\cup V_k'$ is $(\eps,1/3,0)$-super-regular in red.

Take $\delta = \frac{1}{2\Delta}$. We classify the vertices of $V'$
based on their red degrees. We say that a vertex $v$ is {\em good
for} $i$ if either $v \in V_i'$ or the red degree of $v$ to $V_i'$
is at least $\delta|V_i'|/2$, and we say that it is {\em good} if it
is good for every $V_i'$. By regularity, at most an $\eps$ fraction
of the vertices of $V_j'$ fail to be good for $i$, and so at most an
$\eps$ fraction of the vertices of $V'$ fail to be good for $V_i'$.
Define, for each $i$, $B_i$ to be the set of vertices that are not
good for $i$ but are good for every $j$ that is smaller than $i$. By
construction, the $B_i$'s partition the vertices of $V'$ that are
not good, and every vertex in $B_i$ has red degree to $V_i'$ at most
$\delta|V_i'|/2$. Remove the vertices in $\cup_{i=1}^k B_i$ from the
cylinder. Denote the resulting partite sets by $V_i'', 1\leq i \leq
k$. Since
$$|V_i''| \geq (1-\epsilon k) |V_i'| \geq (1-\delta/2)|V_i'|,$$
we have that every vertex in $B_i$ has red degree at most
$\delta|V_i''|$ to $V_i''$. Therefore, it has blue degree at least
$(1-\delta)|V_i''|$ to $V_i''$. Furthermore, since we removed at
most $k\eps |V_i'|$ ($\ll \delta |V_i'|/4$ if $C_2$ is large enough)
vertices from each $V_i'$ and the remaining vertices were all good,
the Slicing Lemma (Lemma \ref{regu}) implies that the resulting
cylinder $V'' = V_1'' \cup \cdots \cup V_k''$ is
$(2\eps,1/4,\delta/4)$-super-regular in red.

Define a $\Delta$-bounded sequence ${\cal{F}}_2'$ of graphs by
taking, for each $m$, a graph from ${\cal{F}}_2$ on $\left\lceil
\frac{\chi_2}{\chi_2-1}m \right\rceil$ vertices, taking an equitable
proper coloring of this graph into $\chi_2$ parts, and removing a
part of size $\lceil m/(\chi_2-1) \rceil$. The resulting graph has
$m$ vertices and an equitable proper coloring into $(\chi_2-1)$
parts. Therefore, by induction, we may partition each $B_i$ into at
most $2^{C(\chi_1 + \chi_2 + \Delta - 1) \log{\Delta}}$ red copies
of graphs from ${\cal{F}}_1$ and blue copies of graphs from
${\cal{F}}_2'$. Denote by $F_1',F_2',\ldots,F_\ell'$ the blue copies
of graphs from ${\cal{F}}_2'$ in such a partition. Each $F_i'$ may
be obtained from some $F_i$ in ${\cal{F}}_2$ by removing an
equitable color class $S_i$ from $F_i$. Define a bipartite graph $H$
whose vertex sets are $A=V(F_1')\cup V(F_2') \cup \cdots \cup
V(F_\ell')$ and $B=S_1 \cup S_2 \cup \cdots \cup S_\ell$ so that the
edges leaving $V(F_i')$ are to $S_i$ and, along with the edges of
the monochromatic $F_i'$, form a copy of $F_i$. Note that $|V(F_i)|
\leq 2|V(F_i')|+1 \leq 3 |V(F_i')|$, and so $|B| \leq 3 |B_i| \leq 3
k \eps |V_i'|$. By Lemma \ref{bipartite}, there is an embedding of
$B$ into some $V_i''' \subseteq V_i''$ so that, along with the
identity embedding on $A$, it forms a homomorphism from $H$ into the
blue edges of $G$. This embedding extends every monochromatic copy
of a graph in ${\cal{F}}_2'$ to a monochromatic copy of a graph in
${\cal{F}}_2$.

The only vertices we have not covered with monochromatic copies of
graphs from ${\cal{F}}_1$ or ${\cal{F}}_2$ are the vertices in each
set of the form $V_i'' \setminus V_i'''$. Since again we removed an
additional at most $k\eps |V_i'|$ ($\ll \delta |V_i''|/8$ if $C_2$
is large enough) vertices from each $V_i''$, the Slicing Lemma
(Lemma \ref{regu}) again implies that the remaining cylinder is
$(2\eps,1/4,\delta/8)$-super-regular in red. Then this remaining
cylinder may be covered by at most $k+1=\Delta+3$ red subgraphs from
${\cal{F}}_1$ by Lemma \ref{cylinder}. Note that the conditions of
the lemma are satisfied if $C_2$ is a sufficiently large absolute
constant. Indeed, we have to check the following inequality
$$\eps = \frac{1}{2^{C_2\Delta}} < \frac{1}{4}\left(\frac{(\frac{1}{32\Delta})
(\frac{1}{8})^{\Delta}}{\Delta(\Delta+2)}\right)^{C_{BL}} =
\frac{1}{4(32\Delta^2 (\Delta + 2) 8^{\Delta})^{C_{BL}}},$$ which is
true if $C_2$ is sufficiently large compared to $C_{BL}$.

The total number of monochromatic subgraphs used in the
partition is at most
$$2^{C_1\Delta
\log{\Delta}}\log(2/(\eps^2 \eta)) + (\Delta+2) 2^{C(\chi_1 + \chi_2
+ \Delta - 1) \log{\Delta}} + (\Delta+3) =$$ $$= 2^{C_1\Delta
\log{\Delta}}\log(2/(\eps^2 \eta)) +
\frac{(\Delta+2)}{2^{C\log{\Delta}}} 2^{C(\chi_1 + \chi_2 + \Delta)
\log{\Delta}} + (\Delta+3) \leq 2^{C(\chi_1 + \chi_2 + \Delta)
\log{\Delta}},$$ as desired, if $C$ is sufficiently large compared
to $C_1$ and $C_2$. Indeed, here by using (\ref{para}) we have
$$ \log(2/(\eps^2 \eta)) = \log(2/(\eps^2)) +  \log(1/\eta) =$$ $$=
(2C_2\Delta + 1) \log{2} + \log{2} + 2(\Delta + 2) \log{2} +
2^{4(\Delta + 2)+5(C_2\Delta+1)}(C_2\Delta + 1) \log{2} \ll
2^{C\Delta}.$$ \qed

\section{Proof of Theorem \ref{bip-tetel}}

The proof is almost identical to the proof of Theorem \ref{tetel}
above. First, a major difference is that in Lemma \ref{ramsey} we
may use the better Ramsey bound for bipartite graphs, thus giving us
the improved bound $2^{C_1\Delta}\log(1/\eps)$. Second, here the
induction has only one step, after one step the chromatic number
goes down to 1 and we may cover each $B_i$ by one graph. This gives
us the bound
$$2^{C_1\Delta}\log(2/(\eps^2 \eta)) + (\Delta + 2) + (\Delta + 3) \leq
2^{C\Delta},$$ as desired. Note that, rather than using the
Hajnal-Szemer\'edi theorem, given any $\Delta$-bounded bipartite
graph sequence ${\cal F}$, we may create a new sequence $G_n$ of
$\Delta$-bounded graphs that have proper equitable $2$-colorings so
that $G_n$ is a union of at most $3$ graphs from ${\cal F}$; we do
this by taking $G_n$ to be two copies of $F_{n/2}$ if $n$ is even
and two copies of $F_{(n-1)/2}$ and a copy of $F_1$ if $n$ is odd.
This would allow us to work with $k=3$ (instead of $\Delta + 2$), so
we may find a regular $3$-cylinder in Step 1. \qed

We should also note the above argument is not particular to bipartite graphs; it works for $\chi$-partite graphs if $\chi$ is a constant (e.g. tripartite graphs); for any constant $\chi$ we get a constant $C(\chi)$ and the bound above becomes $2^{C(\chi)\Delta}$.

\section{Proof of Theorem \ref{lower}}

We wish to show that there exists
a $\Delta$-bounded bipartite sequence ${\cal{F}} =
\{F_1,F_2,\ldots\}$ and, for $n$ sufficiently large, a
two-edge-coloring of $K_n$ that cannot be partitioned into fewer
than $2^{\Omega(\Delta)}$ monochromatic copies of graphs from
$\cal{F}$. To see this, for every $n$ take $G_n$ to be a bipartite
graph on $n$ vertices of degree at most $\Delta$ and, for $n$
sufficiently large, with Ramsey number at least
$2^{\Omega(\Delta)}n$, as given by the result of Graham, R\"odl and
Ruci\'nski \cite{GRR}. We define $F_{2^i}$ recursively; take
$F_{2^0} = G_1$. Then define $F_{2^i}$ to be the disjoint union of
$F_{2^{i-1}}$ with $G_{2^{i-1}}$. For integers of the form $2^i+j$
with $j < 2^i$, define $F_{2^i+j}$ to be the disjoint union of
$F_{2^i}$ with an independent set on $j$ vertices. Under this
definition, each $F_n$ is a bipartite graph on $n$ vertices with
maximum degree at most $\Delta$. Furthermore, for $n_0 < n_1$,
$F_{n_0}$ is a subgraph of $F_{n_1}$. Finally, taking $i$ to be the
largest integer with $2^i \leq n$, $F_n$ contains a copy of
$G_{2^{i-1}}$ and so has Ramsey number at least
$2^{\Omega(\Delta)}2^{i-1} = 2^{\Omega(\Delta)}n$ (for $n$
sufficiently large). Take ${\cal{F}}=\{F_1,F_2,\ldots\}$. Now, for
$N$ sufficiently large, take a $2$-edge-coloring of a complete graph
on $2^{\Omega(\Delta)}N$ vertices without a monochromatic copy of
$F_N$ (this is possible by the condition on the Ramsey number).
Since the sequence of graphs is increasing, this coloring also does
not contain a monochromatic copy of any $F_n$ for $n>N$. Therefore,
any partition of the vertex set into monochromatic copies of graphs
from $\cal{F}$ must use at least $2^{\Omega(\Delta)}$ such copies.
\qed

\section{Concluding Remarks}

There are various interesting potential generalizations of Theorem
\ref{tetel}. One may ask if the theorem holds for $r$ colors for any
positive integer $r$.

\begin{conjecture}\pont
For every positive integer $r$ there exists a constant $C_r$
(depending on $r$) such that, for every $\Delta$-bounded sequence
$\cal{F}$, every $r$-edge-colored complete graph can be partitioned
into at most $2^{\Delta^{C_r}}$ vertex disjoint monochromatic graphs
from $\cal{F}$.
\end{conjecture}

Since bounds on Ramsey numbers were key in proving the theorem for
$r=2$, it is worth noting that Conlon, Fox, and Sudakov \cite{CFS}
proved that, for any fixed number of colors $r$, for any graph $G$
on $n$ vertices of maximum degree $\Delta$ the Ramsey number on $r$
colors $R_r(G)$ is at most $2^{C_r\Delta^2}n$. The primary
difficulty is replacing the step that uses Lemma \ref{bipartite} to
account for the larger number of colors.

Recently B\"ottcher, Kohayakawa, Taraz, and W\"urfl \cite{BKTW}
proved a generalization of the Blow-up Lemma for graphs of bounded
arrangeability without vertices of large degree. An
$a$-arrangeable graph is one in which the vertices may be ordered
such that the neighbors to the right of any vertex $v$ have at most
$a$ neighbors to the left of $v$ in total. They generalize the
Blow-up Lemma from graphs of bounded degree to $n$-vertex graphs of bounded
arrangeability with maximum degree at most $\frac{\sqrt{n}}{\log
n}$. Furthermore, Chen and Schelp \cite{CS} proved that for every
$a$ there is some constant $C(a)$ so that the Ramsey number of any
$a$-arrangeable graph on $n$ vertices is at most $C(a)n$. The best
bound that is known for $C(a)$, again due to Graham, R\"odl and
Ruci\'nski \cite{GRR}, is $C(a) \leq 2^{Ca(\log a)^2}$. One may hope
to combine these two results to get another possible generalization
of Theorem \ref{tetel}.

We say a sequence ${\cal F}=\{F_1,F_2,\ldots\}$ is
$a$-nicely-arrangeable if each $F_n$ is a graph on $n$ vertices that
is $a$-arrangeable with maximum degree at most
$\frac{\sqrt{n}}{\log(n)}$. Using techniques similar to those found in this paper, one can prove:
\begin{theorem}\pont
There exists an absolute constant $C$ so that, for every positive integer $a$ and
every $a$-nicely-arrangeable sequence $\cal{F}$, every
$2$-edge-colored complete graph can be partitioned into at most
$2^{Ca^6}$ vertex disjoint monochromatic graphs from $\cal{F}$.
\end{theorem}

The primary change from the techniques in this paper necessary to prove the above theorem is to adapt the use of Lemma 8. Currently, we use it to take a pair of graph sequences ${\cal F}_1$ and ${\cal F}_2$ along with a nearly-complete bipartite graph, recurse on one of the parts of the bipartite graph to find copies of graphs either from ${\cal F}_1$ or from graphs obtained from ${\cal F}_2$ by removing a color class, and extend the smaller graphs using the nearly-complete bipartite graph. Instead, in the proof for arrangeable graphs, we recurse on one of the parts of the nearly-complete bipartite graph to find copies of graphs from ${\cal F}_1$ and ${\cal F}_2$ (without any color class removed), along with another nearly-complete bipartite graph. This gives a nearly-complete tripartite graph, and we may continue to recurse until we have a nearly-complete multipartite graph into which we may embed our $a$-nicely-arrangeable graphs.

Finally, let us mention that since by now both the Regularity Lemma
and the Blow-up Lemma has been generalized to hypergraphs (see
\cite{RDLS} and \cite{KEE}, respectively), perhaps we can generalize
our result to hypergraphs as well.\\

\noindent {\bf Acknowledgements.} The first author is indebted to Jacob Fox and the second author to
Andr\'as Gy\'arf\'as and Endre Szemer\'edi for helpful conversations on the topic.

\end{document}